\documentclass[twoside]{article}
\title{First, second, and third change of rings theorems for   Gorenstein homological dimensions}
\date{}
\usepackage[all]{xy}
\usepackage{amsfonts}
\usepackage{amsmath}
\usepackage{amssymb}
\usepackage{latexsym}
\newtheorem{thm}{\bf Theorem}[section]

\newtheorem{lem}[thm]{\bf Lemma}
\newtheorem{prop}[thm]{\bf Proposition}
\newtheorem{defn}[thm]{\bf Definition}

\catcode`\ç=13
\defç{\c{c}}
\catcode`\é=13
\defé{\'e}
\catcode`\à=13
\defà{\`a}
\catcode`\è=13
\defè{\`e}
\catcode`\â=13
\defâ{\^a}
\catcode`\ù=13
\defù{\`u}
\catcode`\ê=13
\defê{\^e}
\catcode`\î=13
\defî{\^\i}
\catcode`\ô=13
\defô{\^o}
\newcommand{\field}[1]{\mathbb{#1}}

\newcommand{\Z }{\field{Z}}

\def\proof{{\parindent0pt {\bf Proof.\ }}}

\def\pd{{\rm pd}}
\def\id{{\rm id}}

\def\Gpd{{\rm Gpd}}
\def\Gid{{\rm Gid}}

\def\Im{{\rm Im}}

\def\Ker{{\rm Ker}}

\def\Ext{{\rm Ext}}
\def\Tor{{\rm Tor}}
\def\Hom{{\rm Hom}}

\def\sup{{\rm sup}}

\newcommand{\ZZ}{{\cal Z}}
\newcommand{\rb}{{ \bar{R}}}

\newcommand{\cqfd}
{\hspace{1cm}
\rule{2mm}{2mm}%
\medbreak%
\par%
}
\begin{document}
\thispagestyle{empty}
%%%%%%%%%%%%%%%%%%%%%%%%%%%%%%%%%%%%%%%%%%%%%%%%%%%%%%%%%
%%%%%%%%%%%%%%%%%%%%%%%%%%%%%%%%%%%%%%%%%%%%%%%%%%%%%%%%%
%%%%%%%%%%%%%%%%%%%%%%%%%%%%%%%%%%%%%%%%%%%%%%%%%%%%%%%%%
%%%TITLE%%%%%%%%%%%%%%%%%%%%%%%%%%%%%%%%%%%%%%%%%%%%%%%%%
\maketitle \vspace*{-2cm}
\begin{center}{\large\bf Driss Bennis and Najib Mahdou}
%%%%%%%%%%%%%%%%%%%%%%%%%%%%%%%%%%%%%%%%%%%%%%%%%%%%%%%%%
%%%%%%%%%%%%%%%%%%%%%%%%%%%%%%%%%%%%%%%%%%%%%%%%%%%%%%%%%
%%%%%%%%%%%%%%%%%%%%%%%%%%%%%%%%%%%%%%%%%%%%%%%%%%%%%%%%%
%%%NAMES%%%%%%%%%%%%%%%%%%%%%%%%%%%%%%%%%%%%%%%%%%%%%%%%%

\bigskip
%%%%%%%%%%%%%%%%%%%%%%%%%%%%%%%%%%%%%%%%%%%%%%%%%%%%%%%%%
%%%%%%%%%%%%%%%%%%%%%%%%%%%%%%%%%%%%%%%%%%%%%%%%%%%%%%%%%
%%%%%%%%%%%%%%%%%%%%%%%%%%%%%%%%%%%%%%%%%%%%%%%%%%%%%%%%%
%%%%%%%%%%%%ADDRESSES%%%%%%%%%%%%%%%%%%%%%%%%%%%%%%%%%%%%%%%%%%%%%
\small{Department of Mathematics, Faculty of Science and
Technology of Fez,\\ Box 2202, University S. M.
Ben Abdellah Fez, Morocco, \\[0.12cm] driss\_bennis@hotmail.com\\
mahdou@hotmail.com}
\end{center}

\bigskip\bigskip
%%%%%%%%%%%%%%%%%%%%%%%%%%%%%%%%%%%%%%%%%%%%%%%%%%%%%%%%%
%%%%%%%%%%%%%%%%%%%%%%%%%%%%%%%%%%%%%%%%%%%%%%%%%%%%%%%%%
%%%%%%%%%%%%%%%%%%%%%%%%%%%%%%%%%%%%%%%%%%%%%%%%%%%%%%%%%
%%%ABSTRACT%%%%%%%%%%%%%%%%%%%%%%%%%%%%%%%%%%%%%%%%%%%%%%
\noindent{\large\bf Abstract.}   In this paper, we investigate the
change of rings theorems for the Gorenstein dimensions over
arbitrary rings. Namely, by the use of the notion of strongly
Gorenstein modules, we extend the well-known first, second, and
third change of rings theorems for the classical projective and
injective dimensions to the Gorenstein projective and injective
dimensions, respectively. Each of the results established in this
paper for the Gorenstein projective dimension  is a generalization
of a $G$-dimension of a finitely generated module $M$ over a
noetherian ring $R$.\bigskip

%%%%%%%%%%%%%%%%%%%%%%%%%%%%%%%%%%%%%%%%%%%%%%%%%%%%%%%%%
\medskip


\bigskip\bigskip

%%%%%%%%%%%%%%%%%%%%%%%%%%%%%%%%%%%%%%%%%%%%%%%%%%%%%%%%%
%%%INTRODUCTION%%%%%%%%%%%%%%%%%%%%%%%%%%%%%%%%%%%%%%%%%%
\begin{section}{Introduction}  All rings considered in this paper
are assumed to be commutative with a unit; in particular, $R$
denotes such a ring. All modules are assumed to be
unitary.\bigskip

\textbf{Setup and Notations.} Let $M$  be an $R$-module. An
element $x$ of $R$ is said to be $M$-regular, if $x \not \in
\ZZ_R(M)=\{r\in R\,|\, \exists\; m\in M-\{0\}, \ rm=0\}$ the set
of zero-divisors on $M$. A sequence $x_1,...,x_n$ of elements in
$R$ is called an $M$-sequence, if $(x_1,...,x_n)M\not = M$ and
$x_m \not \in \ZZ_R(M/(x_1,...,x_{m-1})M)$ for $m=1,...,n$ (see
for example \cite{Kap} and \cite{LW}). For an element $x$ of $R$
which is neither a zero-divisor nor a unit, we use $\rb$ to denote
the
quotient ring $R/xR$.\\
We use $\pd_R(M)$ and $\id_R(M)$ to denote, respectively, the
classical projective  and injective
dimension of $M$.\\
We assume that the reader is familiar with the Gorenstein
homological algebra (see references for a background. Namely,
\cite{LW, CFH, CFoxH, Rel-hom, HH}).\bigskip

The Gorenstein homological dimensions theory originated in the
works of Auslander and Bridger  \cite{A1} and \cite{A2}, where
they introduced the G-dimension of any finitely generated module
$M$ and over any Noetherian ring $R$. The G-dimension is analogous
to the G-dimension of a finitely generated module $M$ over a
noetherian ring $R$ and shares some of its principal properties
(see \cite{LW} for more details). However, to complete the
analogy, an extension of the G-dimension to non-necessarily
finitely generated modules is needed. This is done in the 1990's
by Enochs and Jenda \cite{GoIn, GoInPj}  when they defined the
Gorenstein projective dimension, as an extension of the
G-dimension to modules that are not necessarily finitely generated
and over arbitrary  associative rings, and the Gorenstein
injective dimension   as a dual notion of the Gorenstein
projective dimension:

\begin{defn}\label{def1}\begin{itemize}
    \item    \textnormal{An  $R$-module  $M$ is called \textit{Gorenstein projective}
  if there exists an exact sequence of
projective  $R$-modules,
$$\mathbf{P}=\ \cdots\rightarrow P_1\rightarrow P_0 \rightarrow
P_{ -1 }\rightarrow P_{-2 } \rightarrow\cdots,$$ such that  $M
\cong \Im(P_0 \rightarrow P_{ -1 })$ and such that $\Hom_R ( -, Q)
$ leaves the sequence $\mathbf{P}$ exact whenever $Q$ is a
projective  $R$-module.\\
  For a positive integer $n$, we
say that $M$ has \textit{Gorenstein projective dimension} at most
$n$, and we write $\Gpd_R(M)\leq n$ (or simply $\Gpd(M)\leq n$),
if there is an exact sequence of  $R$-modules,
$$ 0 \rightarrow G_n\rightarrow \cdots \rightarrow G_0\rightarrow
M \rightarrow 0,$$ where each $G_i$ is Gorenstein projective.}
\item  \textnormal{Dually, the \textit{Gorenstein injective} module
 is defined, and so the \textit{Gorenstein
injective dimension}, $\Gid_R(M)\leq n$, of an $R$-module $M$ is
defined.}\end{itemize}
\end{defn}

The Gorenstein homological dimensions have been extensively
studied by many others, who proved that these dimensions share
many nice properties of the classical homological dimensions (see,
for instance, \cite{LW,CFoxH, Rel-hom}). In particular, several
classical change of rings theorems have been extended to the
Gorenstein homological dimensions setting (see, for example,
\cite{LW}, \cite{CHH}, and  \cite{KY2}). However, most of those
results are proved over Noetherian local rings, and in many cases
the rings are in addition homomorphic images of Gorenstein local
rings. In this paper, we generalize some change of rings results
for Gorenstein homological dimensions by relaxing the conditions
on the underlying rings. We are mainly interested in the
Gorenstein counterpart of the  classical first, second, and third
change of rings theorems (see \cite[Sections 4.4 and 4.5]{Kap} and
\cite[Section 4.3]{Wei}\,\footnote{\;In this paper, we refer to
Weibel's book \cite{Wei} for the projective case and to
Kaplansky's book \cite{Kap} for the injective case.}). Namely, we
extend these classical theorems to the Gorenstein setting.\bigskip

Before investigating the first one, we establish, in Section 3,
the second change of rings theorems for the Gorenstein projective
and
  injective dimensions. Those are Theorems \ref{thm-second-Gp} and
\ref{thm-second-Gi}:\bigskip

\noindent\textbf{Second change of rings theorems for the
Gorenstein projective and injective dimensions.}\\ \textit{Let $M$
be a non-zero $R$-module and let $x\in R$ be both $R$-regular and
$M$-regular. Then,
\begin{enumerate}
    \item  $\Gpd_{\rb}(M/xM)\leq \Gpd_{R}(M).$
    \item  $\Gid_{\rb}(M/xM)\leq \Gid_{R}(M)-1 $, except when $M$
    is a Gorenstein injective $R$-module.
\end{enumerate}}\bigskip

These inequalities are extensions of the classical second change
of rings theorems for the projective and injective dimensions
(\cite[Theorem 4.3.5]{Wei} and \cite[Theorem 205]{Kap}),
respectively. The first one is a generalization of
\cite[Proposition 1.3.6]{LW}, and the second one generalizes
\cite[Corollary 2.3]{KY}  in which $R$ is assumed to be Noetherian
local and $M$ to be finitely generated.\bigskip

Section 4 is devoted to the first change of rings theorems the
Gorenstein projective and injective dimensions. Recall the first
change of rings theorem for the projective and injective
dimensions (\cite[Theorem 4.3.3]{Wei} and \cite[Theorem
202]{Kap}):\bigskip

\noindent\textbf{First change of rings theorems for the projective
and injective dimensions.}\\ \textit{Let $M$ be a non-zero
$R$-module and let $x=x_1,...,x_t$ be an $R$-sequence of elements
in the annihilator of $M$.  Then,
\begin{enumerate}
    \item If  $\pd_{R/(x)}(M)<\infty$, then  $\pd_{R}(M)=\pd_{R/(x)}(M)+t.$
    \item If  $\id_{R/(x)}(M)<\infty$, then  $\id_{R}(M)=\id_{R/(x)}(M)+t.$\\
\end{enumerate}}

The first equality was extended for finitely generated modules and
over  Noetherian local rings to the Gorenstein projective
dimension, as follows \cite[Proposition 1.5.3]{LW}:\bigskip

\noindent\textit{Let $R$ be a  Noetherian local ring. Let $M$ be a
non-zero $R$-module and let $x=x_1,...,x_t$ be an $R$-sequence of
elements in the annihilator
of $M$.  Then,  \\
     If  $\Gpd_{R/(x)}(M)<\infty$, then  $\Gpd_{R}(M)=\Gpd_{R/(x)}(M)+t.$
     \bigskip
 }

It is known that the first change of rings theorem for the
projective dimension does not hold if we remove the condition
$\pd_{R/(x)}(M)<\infty$. For example, we have $\pd_{\Z/4\Z}(
\Z/2\Z)=\infty$ but $\pd_{\Z}( \Z/2\Z)=1$ \cite[Example
4.3.1]{Wei}. In this example, $\Z/2\Z\cong 2\Z/4\Z$ is a
Gorenstein projective $\Z/4\Z$-module. This implies that
$$\Gpd_{\Z}( \Z/2\Z)= \Gpd_{\Z/4\Z}( \Z/2\Z)+1.$$ Similar example
leads Christensen to ask whether the condition
$\Gpd_{R/(x)}(M)<\infty $ in \cite[Proposition 1.5.3]{LW} is
necessary \cite[Remarks 1.5.4]{LW}. Later, in \cite[Theorem
2.2.8]{LW}, he proved that certainly the equality
$\Gpd_{R}(M)=\Gpd_{R/(x)}(M)+t$ holds without
 assuming the condition $\Gpd_{R/(x)}(M)<\infty $. In section
4, we generalize \cite[Theorem 2.2.8]{LW}. Namely, we show that
the equality $\Gpd_{R}(M)=\Gpd_{R/(x)}(M)+t$ holds over arbitrary
rings and for non-necessarily finitely generated modules. Also, we
establish its dual result for the Gorenstein injective dimension.
Those are Theorems \ref{thm-first-GP} and
\ref{thm-first-Gi}:\bigskip

\noindent\textbf{First change of rings theorem for the Gorenstein
projective and injective dimensions.}\\ \textit{Let $M$ be a
non-zero $R$-module and let $x=x_1,...,x_t$ be an $R$-sequence of
elements in the annihilator of $M$.  Then,
   $$\Gpd_{R}(M)=\Gpd_{R/(x)}(M)+t\quad  and\quad
    \Gid_{R}(M)=\Gid_{R/(x)}(M)+t. $$
    }

Finally, in Section 5, we discuss the third change of rings
theorem.\bigskip

In the next section, we give some definitions and results needed
in the rest of the paper.

\end{section}
%%%%%%%%%%%%%%%%%%%%%%%%%%%%%%%%%%%%%%%%%%%%%%%%%%%%%%%%%%%%%%%%%%%%%%%%%%%
%%%%%%%%%%%%%%%%%%%%%%%%%%%%%%%%%%%%%%%%
%%%%%%%%%%%%%%%%%%%%%%%%%%%%%%%%
%%%%%%%%%%%%%%%%%%%%%%%%                  $$$$                      $$$$
%%%%%%%%%%%%%%%%%%%                          Section 2:  Preliminaries
%%%%%%%%%%%%%%%%%%%%%%%%                  $$$$                      $$$$
%%%%%%%%%%%%%%%%%%%%%%%%%%%%%%%%
%%%%%%%%%%%%%%%%%%%%%%%%%%%%%%%%%%%%%%%%
%%%%%%%%%%%%%%%%%%%%%%%%%%%%%%%%%%%%%%%%%%%%%%%%%%%%%%%%%%%%%%%%%%%%%%%%%%%%
%%%%%%%%%%%%%%%%%%%%%%%%%%%%%%%%%%%%%%%%%%%%%%%%%%%%%%%%%%%%
\begin{section}{Preliminaries}
The proofs of almost all theorems given in the present paper are
mainly based on the notion of ``strongly Gorenstein projective and
injective modules''. These modules are  introduced in \cite{BM1},
as particular cases of the Gorenstein projective and injective
modules, respectively, as follows:

\begin{defn}[\cite{BM1}]\label{def-G-str-proj-mod}\begin{enumerate}
    \item \textnormal{An $R$-module $M$ is said to be strongly Gorenstein projective, if
there exists an exact sequence of the form
$$ \mathbf{P}=\
\cdots\stackrel{f}{\longrightarrow}P\stackrel{f}{\longrightarrow}P\stackrel{f}{\longrightarrow}P
\stackrel{f}{\longrightarrow}\cdots,$$
     where $P$ is a projective $R$-module and $f$ is an endomorphism of $P$,
     such that  $M \cong \Im(f)$ and such that $\Hom_R( -, Q) $ leaves the
sequence $\mathbf{P}$ exact whenever $Q$ is a projective
$R$-module.}
\item  \textnormal{Dually the strongly Gorenstein injective $R$-modules is defined.}
\end{enumerate}
\end{defn}

These particular cases of the Gorenstein projective and injective
modules have simpler characterizations:

\begin{prop}[\cite{BM1}, Proposition 2.9 and Remark 2.10
(2)]\label{por-cara-G-str-mod} \begin{enumerate}     \item An
$R$-module $M$ is strongly Gorenstein projective if  and only
    if
there exists a short exact sequence of $R$-modules $$0\rightarrow
M\rightarrow P\rightarrow M\rightarrow 0,$$ where $P$ is a
projective  $R$-module, and  $\Ext^1_R(M,Q)=0$ for any  $R$-module
$Q$ with finite projective dimension (or for any projective
$R$-module $Q$).
    \item An $R$-module $M$ is strongly Gorenstein injective if and only
    if
there exists a short exact sequence of $R$-modules $$0\rightarrow
M\rightarrow I\rightarrow M\rightarrow 0,$$ where $I$ is an
injective $R$-module, and  $\Ext^1_R(E, M)=0$ for any  $R$-module
$E$ with finite injective dimension (or for any injective
$R$-module $E$).
\end{enumerate}
\end{prop}

The principal role of the strongly Gorenstein projective and
injective modules is to give a simple characterization of the
Gorenstein projective and injective modules, respectively, as
follows:

\begin{thm}[\cite{BM1}, Theorem 2.7] \label{thm-car-G-SG}
An $R$-module is Gorenstein projective  (resp., injective) if and
only if  it is a direct summand of a strongly Gorenstein
projective (resp., injective) $R$-module.
\end{thm}

We also need  the following extensions of the well-known standard
(in)equalities for the projective dimension  (\cite[Corollary 2,
p. 135]{Bou}) to the Gorenstein projective dimension.

\begin{lem}\label{lem-standard-ineq-Gp}
Let $0\rightarrow A \rightarrow  B\rightarrow C\rightarrow0$ be a
short exact sequence of $R$-modules. Then,
\begin{enumerate}
    \item  $\Gpd(A)\leq \sup\{\Gpd(B),\Gpd(C)-1\}$ with equality if
    $\Gpd(B)\not= \Gpd(C)$.
    \item $\Gpd(B)\leq \sup\{\Gpd(A),\Gpd(C)\}$ with equality if
    $\Gpd(C)\not= \Gpd(A)+1$.
    \item  $\Gpd(C)\leq \sup\{\Gpd(B),\Gpd(A)+1\}$ with equality if
    $\Gpd(B)\not= \Gpd(A)$.
\end{enumerate}
\end{lem}
\proof Using \cite[Theorems 2.20 and 2.24]{HH} the argument is
analogous to the one of \cite[Corollary 2, p.
135]{Bou}.\cqfd\bigskip

Some of these  (in)equalities are already proved in special cases
(see \cite[Corollary 1.2.9]{LW} and \cite[Proposition
2.18]{HH}).\bigskip

Dually, we can prove the Gorenstein injective version of Lemma
\ref{lem-standard-ineq-Gp} above.

\begin{lem}\label{lem-standard-ineq-Gi}
Let $0\rightarrow A \rightarrow  B\rightarrow C\rightarrow0$ be a
short exact sequence of $R$-modules. Then,
\begin{enumerate}
    \item  $\Gid(A)\leq \sup\{\Gid(B),\Gid(C)+1\}$ with equality if
    $\Gid(B)\not= \Gid(C)$.
    \item $\Gid(B)\leq \sup\{\Gid(A),\Gid(C)\}$ with equality if
    $\Gid(A)\not= \Gid(C)+1$.
    \item  $\Gid(C)\leq \sup\{\Gid(B),\Gid(A)-1\}$ with equality if
    $\Gid(B)\not= \Gid(A)$.
\end{enumerate}
\end{lem}

\end{section}

%%%%%%%%%%%%%%%%%%%%%%%%%%%%%%%%%%%%%%%%%%%%%%%%%%%%%%%%%%%%%%%%%%%%%%%%%%%
%%%%%%%%%%%%%%%%%%%%%%%%%%%%%%%%%%%%%%%%
%%%%%%%%%%%%%%%%%%%%%%%%%%%%%%%%
%%%%%%%%%%%%%%%%%%%%%%%%                  $$$$                      $$$$
%%%%%%%%%%%%%%%%%%%                          Section 3:  Second
%%%%%%%%%%%%%%%%%%%%%%%%                  $$$$                      $$$$
%%%%%%%%%%%%%%%%%%%%%%%%%%%%%%%%
%%%%%%%%%%%%%%%%%%%%%%%%%%%%%%%%%%%%%%%%
%%%%%%%%%%%%%%%%%%%%%%%%%%%%%%%%%%%%%%%%%%%%%%%%%%%%%%%%%%%%%%%%%%%%%%%%%%%%
%%%%%%%%%%%%%%%%%%%%%%%%%%%%%%%%%%%%%%%%%%%%%%%%%%%%%%%%%%%%
\begin{section}{Second change of rings theorems for the Gorenstein projective and injective
dimensions} This section is devoted to the second change of rings
theorems for the Gorenstein projective and injective dimensions.
We begin with the Gorenstein projective case, which is an
extension of the well-known ``second change of rings theorem for
the projective dimension'' \cite[Theorem 4.3.5]{Wei} and it is a
generalization of \cite[Proposition 1.3.6]{LW}.

\begin{thm}\label{thm-second-Gp}
Let $M$ be a non-zero $R$-module and let $x\in R$ be both
$R$-regular and $M$-regular. Then,
  $$\Gpd_{\rb}(M/xM)\leq \Gpd_{R}(M).$$
  \end{thm}
\proof Since $\Tor_{i}^{R}(M,\rb)=0$ for all $i>0$ (by
\cite[Examples (1), p. 102]{Bou}), the inequality follows by the
following general change of rings result.\cqfd

\begin{lem}\label{lem-second-Gp}
Let $R\rightarrow S$ be a ring homomorphism with $\pd_R(S)<\infty$
and let $M$ be a non-zero $R$-module. If $\Tor_{i}^{R}(M,S)=0$ for
all $i>0$, then $\Gpd_{S}(M\otimes_R S)\leq \Gpd_{R}(M)$.
    \end{lem}
\proof To prove the inequality, we may assume that
$\Gpd_{R}(M)=m<\infty$. First assume that $M$ is a Gorenstein
projective $R$-module. We prove that $M\otimes_R S$ is a
Gorenstein projective $S$-module. For that, we may assume, by
\cite[Theorem 2.5]{HH} and Theorem \ref{thm-car-G-SG}, that $M$ is
strongly Gorenstein projective. Then, from Proposition
\ref{por-cara-G-str-mod}, there exists a short exact sequence of
$R$-modules $$0\rightarrow M\rightarrow P\rightarrow M\rightarrow
0,$$ where $P$ is projective, and $\Ext^1_R(M,Q)=0$ for any
$R$-module $Q$ with finite projective dimension. By the short
exact sequence above, $\Tor_{i}^{R}(M,S)\cong\Tor_{i+1}^{R}(M,S)$
for all $i>0$. Then, $\Tor_{i}^{R}(M,S)=0$  for all $i>0$ (since
$\pd_R(S)<\infty$), and so we get the following exact sequence of
$S$-modules
$$0\rightarrow M\otimes_R S\rightarrow P\otimes_R S\rightarrow
M\otimes_R S\rightarrow 0.$$ On the other hand, consider a
projective $S$-module $Q$. Then, $\pd_R(Q)$ is finite (by
\cite[Exercise 5, p. 360]{CE}), and so, by \cite[Proposition
4.1.3]{CE} and the fact that $M$ is a Gorenstein projective
$R$-module
$$\Ext^1_S(M\otimes_R S, Q)\cong \Ext^1_R(M , Q)=0.$$
Therefore, $M\otimes_R S$ is a strongly Gorenstein projective
$S$-module (by Proposition
\ref{por-cara-G-str-mod}).\\
Assume now that $m=\Gpd_{R}(M)>0 $. Then, there exists a short
exact sequence of $R$-modules,  $$0\rightarrow K\rightarrow
F\rightarrow M\rightarrow 0,$$ where $F$ is free and, by Lemma
\ref{lem-standard-ineq-Gp},  $\Gpd_R(K)= m-1$. This, implies, by
the induction hypothesis, that $\Gpd_{S}(K\otimes_R S)\leq m-1$.
On the other hand, since $\Tor_1^R(M,S)=0$, we get a short exact
sequence of $S$-modules
$$0\rightarrow K\otimes_R S\rightarrow F\otimes_R S\rightarrow
M\otimes_R S\rightarrow 0.$$ Therefore, $\Gpd_{S}(M\otimes_R
S)=\Gpd_R(K\otimes_R S)+1 \leq m $ (by Lemma
\ref{lem-standard-ineq-Gp}, since $F\otimes_R S$ is a free
$S$-module), which completes the proof.\cqfd\bigskip

Dually we prove the following  general change of rings result,
which is a generalization of \cite[Theorem 204]{Kap}.

\begin{lem}\label{lem-second-Gi}
Let $R\rightarrow S$ be a ring homomorphism with $\pd_R(S)<\infty$
and let $M$ be a non-zero $R$-module. If $\Ext_{R}^{i}(S,M)=0$ for
all $i>0$, then $ \Gid_{S}(\Hom_R(S , M))\leq \Gid_{R}(M)$.
\end{lem}

Now, we give an extension of the ``second change of rings theorem
for the injective dimension'' \cite[Theorem 205]{Kap}. This
theorem generalizes \cite[Corollary 2.3]{KY}, which is obtained
for finitely generated modules and over Noetherian local rings.

\begin{thm}\label{thm-second-Gi}
Let $M$ be a non-zero $R$-module and let $x\in R$ be both
$R$-regular and $M$-regular. Then,  $$ \Gid_{\rb}(M/xM)\leq
\Gid_{R}(M)-1,$$ except when $M$ is a Gorenstein injective
$R$-module (in which case $xM=M$).
\end{thm}
\proof  First note that if $M$ is a Gorenstein injective
$R$-module, then it is divisible; i.e., $xM=M$ for every
non-zero-divisor element $x$ of
$R$ (since it is a quotient of an injective module).\\
Now, to prove the inequality, we may assume that $\Gid_{R}(M)=n$
with $1\leq n < \infty$. Then, there exists a short  exact
sequence of $R$-modules,
$$0\rightarrow M\rightarrow E\rightarrow I\rightarrow 0,$$ where $E$
is injective and, by Lemma \ref{lem-standard-ineq-Gi}, $\Gid_R
(I)=n-1$. From \cite[Examples (1), p. 102]{Bou}, we have that
$\Tor_1^R(N,\rb)=N_x$ for  all $R$-modules $N$, where $N_x$
denotes the submodule of $N$ annihilated by $x$. Then, Tensorising
the  short sequence above by $\rb$ we get the following short
exact sequence of $\rb$-modules
$$0\rightarrow E_x\rightarrow I_x\rightarrow M/xM\rightarrow  0.$$
By \cite[Theorem 204]{Kap}, $E_x$ is an injective $\rb$-module,
then the short exact sequence splits and so $M/xM$ is a direct
summand of $I_x$. Thus, $\Gid_{\rb}(M/xM)\leq \Gid_{\rb}( I_x)$
(by Lemma \ref{lem-standard-ineq-Gi}). Finally, the desired
inequality $\Gid_{\rb}(M/xM)\leq n-1$ is obtained by the
isomorphism $\Hom_R(\rb, I)\cong I_x$ and by Lemma
\ref{lem-second-Gi} (In fact,  from \cite[Examples (1), p.
102]{Bou} and since $I$ is divisible as a quotient of the
injective (then divisible) $R$-module $E$, $I$ satisfies the
hypothesis $\Ext_{R}^{i}(\rb,I)=0$ for all $i>0$ of Lemma
\ref{lem-second-Gi}).\cqfd

\end{section}
%%%%%%%%%%%%%%%%%%%%%%%%%%%%%%%%%%%%%%%%%%%%%%%%%%%%%%%%%
%%%%%%%%%%%%%%%%%%%%%%%%%%%%%%%%%%%%%%%%%%%%%%%%%%%%%%%%%
%%%%%%%%%%%%%%%%%%%%%%%%%%
%%%%%%%%%%%%%%%%%%%%%        Section 4:   First
%%%%%%%%%%%%%%%%%%%%%
%%%%%%%%%%%%%%%%%%%%%%%%%%%%%%%%%%%%%%%%%%%%%%%%%%%%%%%%%%%%%
%%%%%%%%%%%%%%%%%%%%%%%%%%%%%%%%%%%%%%%%%%%%%%%%%%%%%%%%%%%%%%%
\begin{section}{First change of rings theorems for the Gorenstein projective and injective
dimensions} In this section, we investigate the first change of
rings theorems for the Gorenstein projective and injective
dimensions.
We begin with the Gorenstein projective case.\\
\indent As mentioned in the introduction, the classical first
change of rings theorem for the projective dimension was already
generalized and extended to the Gorenstein projective dimension
for finitely generated modules and over Noetherian rings
\cite[Theorem 2.2.8]{LW}. Here, we show, by completely different
arguments, that \cite[Theorem 2.2.8]{LW} can  be generalized to
arbitrary rings and to non-necessarily finitely generated modules.
That is the following first change of rings theorem for the
Gorenstein projective dimension.

\begin{thm} \label{thm-first-GP}
Let $M$ be a non-zero $R$-module and let $x=x_1,...,x_t$ be an
$R$-sequence of elements in the annihilator of $M$. Then
 $$\Gpd_{R}(M)=\Gpd_{R/(x)}(M)+t.$$
 In particular, $\Gpd_{R}(M)$ and $\Gpd_{R/(x)}(M)$ are
simultaneously finite.
\end{thm}
\proof By induction on the length $t$ of the $R$-sequence $x$, it
suffices to prove the case $t=1$, such that we write $x_1=x$ and
$\rb=R/xR$.\\
First we prove that $\Gpd_{\rb}(M)<\infty$ implies
$\Gpd_{R}(M)=\Gpd_{\rb}(M)+1.$ This gives the first inequality
$\Gpd_{R}(M)\leq \Gpd_{\rb}(M)+1.$\\
Assume then that $\Gpd_{\rb}(M)=n<\infty$. We prove the equality
$\Gpd_{R}(M)=n+1$ by induction on $n$. Assume first that $M$ is a
Gorenstein projective $\rb$-module and prove that $\Gpd_{R}(M)=1$.
Note that $M$ can not be a Gorenstein projective $R$-module.
Indeed, any Gorenstein projective $R$-module can be embedded in a
free $R$-module and therefore can not have the $R$-regular element
$x$ as a zero-divisor. So, it suffices to prove $\Gpd_{R}(M)\leq
1$. By Theorem \ref{thm-car-G-SG} and \cite[Proposition 2.19]{HH},
we may consider $M$ to be a strongly Gorenstein projective
$\rb$-module. Thus, from Proposition \ref{por-cara-G-str-mod},
there exists a short exact sequence of $\rb$-modules
$$0\rightarrow M\rightarrow P\rightarrow M\rightarrow 0,$$ where $P$
is a projective $\rb$-module, and $\Ext^1_{\rb}(M,Q)=0$ for any
projective $\rb$-module $Q$. From \cite[Lemma 6.20]{Rot}, there is
a commutative diagram:
$$\begin{array}{ccccccccc}
    &  & 0 &  & 0&  &0 &  &  \\
    &  & \downarrow &  &\downarrow&  &\downarrow &  &  \\
   0 & \rightarrow & P_1 &\rightarrow  &P'& \rightarrow &P_1 &\rightarrow  & 0 \\
    &  & \downarrow &  &\downarrow&  &\downarrow &  &  \\
   0 & \rightarrow & P_0&\rightarrow  &P_0\oplus P_0& \rightarrow &P_0 &  \rightarrow&  0\\
    &  & \downarrow &  &\downarrow&  &\downarrow &  &  \\
   0 & \rightarrow &M &\rightarrow  &P& \rightarrow &M & \rightarrow & 0  \\
    &  & \downarrow &  &\downarrow&  &\downarrow &  &  \\
    &  & 0 &  & 0&  &0 &  &
\end{array}$$
where $P_0$ is a projective $R$-module. Then, $P'$ is  a
projective $R$-module (by \cite[Theorem 4.3.3]{Wei}). On the other
hand, by the first vertical sequence in the diagram above, we have
the isomorphism
$$\Ext^1_R(P_1,Q)\cong \Ext_{R}^{2}(M,Q)$$  for every $R$-module
$Q$. If such an $R$-module $Q$ is projective (and so $Q/xQ$ is a
projective $\rb$-module), then the Rees's theorem \cite[Theorem
9.37]{Rot} gives
$$\Ext_{R}^{2}(M,Q)\cong\Ext^1_{\rb}(M,Q/xQ)=0.$$ Then,
$\Ext^1_R(P_1,Q)=0$ and so $P_1$ is a strongly Gorenstein
projective $R$-module (by Proposition \ref{por-cara-G-str-mod}).
This means, by the diagram above, that $\Gpd_{R}(M)\leq 1$. Thus,
the proof of the case $n=0$ is finished.\\[0.2cm]
Now, assume that $\Gpd_{\rb}(M)=n>0$. Then, there exists a short
exact sequence of $\rb$-modules, $$0\rightarrow K \rightarrow F
\rightarrow M \rightarrow 0,$$  where $F$ is  free and, by Lemma
\ref{lem-standard-ineq-Gp}, $\Gpd_{\rb}(K)=n-1$. Hence, by
induction, $\Gpd_{R}(K)=n$. Therefore, $\Gpd_R(M)= \Gpd_{R}(K)+1
=n+1$ (by Lemma \ref{lem-standard-ineq-Gp} and since
$\Gpd_{R}(F)=\pd_{R}(F)=1$ (by \cite[Theorem 4.3.3]{Wei} and
\cite[Proposition 2.27]{HH})).
This completes the proof of the first part.\\[0.2cm]
Now, it remains to prove the inequality $\Gpd_{\rb}(M)+1\leq
\Gpd_{R}(M) $. For that, we may assume that
$\Gpd_{R}(M)=n<\infty$. As shown  in the first part of this proof,
$M$ can not be a Gorenstein projective $R$-module, and so we begin
with the case $ \Gpd_{R}(M)=1$ and we prove that $M$ is a
Gorenstein projective $\rb$-module. This is equivalent, by
\cite[Proposition 2.3]{HH}, to show the following two statements:
\begin{enumerate}
    \item $\Ext^{i}_{\rb}(M, Q)=0$ for every
$i\geq 1$ and every projective $\rb$-module  $Q$.
    \item There exists an exact sequence of $\rb$-modules
    $$\mathbf{P}=\quad 0\rightarrow M \rightarrow  P^0
 \rightarrow P^1\rightarrow\cdots$$  such that $\Hom_R ( -, Q) $ leaves the
sequence $P$ exact whenever $Q$ is a projective $\rb$-module.
\end{enumerate}
To prove the first assertion, it suffices to consider $Q$ to be a
free $\rb$-module. In this case $Q$ is of the form $F/xF$, where
$F$ is a free $R$-module. Therefore, from Rees's theorem
\cite[Theorem 9.37]{Rot} and since  $ \Gpd_{R}(M)=1$, we get for
every $i\geq 1$
$$(\ast)\quad \Ext_{\rb}^{i}(M,F/xF) \cong \Ext_{R}^{i+1}(M,F)=0.$$
Now, we prove the existence of the resolution $\mathbf{P}$. For
that, it suffices to prove the existence of short exact sequences
$$ 0\rightarrow M_j \rightarrow P^j \rightarrow M_{j+1} \rightarrow 0
\qquad j\geq0,$$ where $M_0=M$ and each $P^j$ is a projective
$\rb$-module, such that  $\Ext^1_{\rb}(M_j,Q)=0$
for any $j\geq 0$ and any projective $\rb$-module $Q$.\\
Since $ \Gpd_{R}(M)=1$, there exists a short exact sequence of
$R$-modules
$$0\rightarrow G_1\rightarrow G_0\rightarrow M\rightarrow 0,$$
where $G_0$  is free and   $G_1$ is Gorenstein projective.
Tensorising this sequence by $\rb$, we get the following exact
sequence of $\rb$-module:
$$(\Tor_1^R(\rb,G_0)=)\, 0\rightarrow M \rightarrow G_1/xG_1\rightarrow G_0/xG_0\rightarrow M \rightarrow 0.$$
Setting $N=\Ker( G_0/xG_0\rightarrow M )$  we get two short exact
sequences of $\rb$-modules:
$$\begin{array}{c}
   0\rightarrow M \rightarrow G_1/xG_1\rightarrow N \rightarrow 0\qquad
   \mathrm{and}\qquad
  0\rightarrow N\rightarrow G_0/xG_0\rightarrow M \rightarrow 0.
  \end{array}$$
Since $G_1$ is a Gorenstein projective $R$-module, $x$ is
$G_1$-regular. Then, by Theorem \ref{thm-second-Gp}, $G_1/xG_1$ is
a Gorenstein projective $\rb$-module. Thus, there exists, by
definition, a short exact sequence of $\rb$-modules $$
0\rightarrow
 G_1/xG_1 \rightarrow P^0 \rightarrow H_0 \rightarrow 0,$$
 where $P^0$ is projective and $H_0$ is Gorenstein projective.
 Consider the following pushout diagram:
  $$ \xymatrix{
    &  & 0  \ar[d]        &    0\ar[d]          &    \\
   0 \ar[r] &  M \ar@{=}[d] \ar[r]  & G_1/xG_1 \ar[d] \ar[r]          &   N   \ar@{-->}[d] \ar[r] & 0\\
     0 \ar[r]   & M     \ar[r]  &  P^0      \ar[d] \ar@{-->}[r]    & M_1    \ar[d] \ar[r]      & 0\\
   &   &         H_0  \ar[d]  \ar@{=}[r]& H_0 \ar[d]           &   &  \\
  &  & 0           & 0        &     }$$
Since $G_0$ is a free $R$-module, $G_0/xG_0$ is a free
$\rb$-module, and so $\pd_R(G_0/xG_0)=1$. Then, $\Gpd_R(N)\leq 1$
(by Lemma \ref{lem-standard-ineq-Gp}  and by the short exact
sequence of $\rb$-modules   $0\rightarrow N\rightarrow
G_0/xG_0\rightarrow M \rightarrow 0 $). Then, using the same
argument as above, we get  that   $\Ext^{i}_{\rb}(N, Q)=0$ for
every $i\geq 1$ and every projective $\rb$-module  $Q$.  Thus, the
middle horizontal sequence of the diagram above, $$ 0 \rightarrow
M_0 (=M)\rightarrow P^0 \rightarrow M_1 \rightarrow
0,$$ is  the desired first short exact sequence.\\
Now using the short exact sequence $ 0\rightarrow N\rightarrow
G_0/xG_0\rightarrow M \rightarrow 0$ and the right vertical
sequence in the diagram above, we get the following pushout
diagram:
 $$ \xymatrix{
   &   0      \ar[d]        &  0   \ar[d]           &   &      \\
   0  \ar[r] & N  \ar[d] \ar[r] & G_0/xG_0 \ar@{-->}[d] \ar[r] & M   \ar@{=}[d] \ar[r] & 0\\
     0   \ar[r] & M_1   \ar[d] \ar@{-->}[r] & O_1     \ar[d] \ar[r] & M   \ar[r] & 0\\
    &  H_0   \ar[d] \ar@{=}[r] & H_0     \ar[d]  &       &     &  \\
    &   0             &  0            &   &         }$$
Since  $O_1$  is a Gorenstein projective $\rb$-module (apply
\cite[Theorem 2.5]{HH} to the middle vertical sequence in the
diagram above),  there exists a short exact sequence of
$\rb$-modules
$$ 0\rightarrow
 O_1 \rightarrow P^1 \rightarrow H_1 \rightarrow 0,$$
 where $P^1$ is projective and $H_1$ is Gorenstein projective.
 With this sequence and the middle horizontal sequence in the diagram above we get the following pushout
 diagram:
  $$ \xymatrix{
      &  & 0      \ar[d]      & 0   \ar[d]         &  \\
   0  \ar[r] & M_1    \ar@{=}[d] \ar[r]     & O_1 \ar[d] \ar[r]     & M \ar@{-->}[d] \ar[r]  & 0\\
     0  \ar[r] & M_1  \ar[r] &  P^1     \ar[d] \ar@{-->}[r]& M_2        \ar[d] \ar[r] &     0\\
       &   &  H_1        \ar[d] \ar@{=}[r] & H_1   \ar[d]            &  \\
    &     &  0             &  0           &         }$$
By the right vertical sequence, $\Ext^1_{\rb}(M_2,Q)=0$ for any
projective $\rb$-module $Q$.  Thus the desired second short exact
sequence $0 \rightarrow M_1\rightarrow  P^1 \rightarrow M_2
\rightarrow 0$ is established.\\
Recursively we construct the remains short exact sequences, and
this completes the proof of the case $n=1$.\\[0.2cm]
Finally, assume that $\Gpd_R(M)=n>1$ and consider a short exact
sequence of $\rb$-modules, $$ 0 \rightarrow K\rightarrow
F\rightarrow M\rightarrow 0,$$ where $F$ is  free. Since $
\pd_R(F)=1$ and by Lemma \ref{lem-standard-ineq-Gp},
$\Gpd_R(K)=\Gpd_R(M)-1=n-1$. Thus, the induction hypothesis gives
$\Gpd_{\rb}(K)\leq\Gpd_R(K)-1=n-2$ and therefore
$$\Gpd_{\rb}(M)=\Gpd_{\rb}(K)+1 \leq n-1.$$ This completes the
proof.\cqfd\bigskip

Dually, we get the first change of rings theorem for the
Gorenstein injective dimension.

\begin{thm} \label{thm-first-Gi}
Let $M$ be a non-zero $R$-module and let $x=x_1,...,x_t$ be an
$R$-sequence of elements in the annihilator of $M$. Then,
 $$\Gid_{R}(M)=\Gid_{R/(x)}(M)+t.$$
In particular, $\Gid_{R}(M)$ and $\Gid_{R/(x)}(M)$ are
simultaneously finite.
\end{thm}

\end{section}

%%%%%%%%%%%%%%%%%%%%%%%%%%%%%%%%%%%%%%%%%%%%%%%%%%%%%%%%%
%%%%%%%%%%%%%%%%%%%%%%%%%%%%%%%%%%%%%%%%%%%%%%%%%%%%%%%%%
%%%%%%%%%%%%%%%%%%%%%%%%%%
%%%%%%%%%%%%%%%%%%%%%        Section 4: Third
%%%%%%%%%%%%%%%%%%%%%
%%%%%%%%%%%%%%%%%%%%%%%%%%%%%%%%%%%%%%%%%%%%%%%%%%%%%%%%%%%%%
%%%%%%%%%%%%%%%%%%%%%%%%%%%%%%%%%%%%%%%%%%%%%%%%%%%%%%%%%%%%%%%
\begin{section}{Third change of rings theorem for Gorenstein projective dimension}
Recall the third change of rings theorem for the classical
projective dimension \cite[Theorem 4.3.12 and Remark p. 104]{Wei}:
Let $M$ be a non-zero $R$-module and let $x\in R$ be both
$R$-regular and $M$-regular. If $R$ is Noetherian and $x$ is in
the Jacobson radical of $R$, then $\pd_{\rb}(M/xM)= \pd_{R}(M).$
This theorem is extended to coherent rings and for finitely
presented modules (see \cite[Theorem 3.1.2]{Glaz}). In the
Gorenstein dimensions theory, we find an extension of the third
change of rings theorem for the projective dimension to the
Gorenstein projective dimension of finitely generated modules and
over  Noetherian local rings (see \cite[Corollary 1.4.6]{LW}).
Here, as the classical case, we extend \cite[Corollary 1.4.6]{LW}
as follows:

\begin{thm}\label{thm-third}
Let $R$ be a coherent ring and let $M$ be a non-zero finitely
presented $R$-module. If $x=x_1,...,x_t$ is an $R$-sequence in the
Jacobson radical of $R$ and an $M$-sequence, then
    $$\Gpd_{R/(x)}(M/(x)M)= \Gpd_{R}(M).$$
    \end{thm}
\proof  Using  \cite[Proposition 10.2.6 $(1)\Leftrightarrow
(10)$]{Rel-hom}, the proof is the same as the one of
\cite[Corollary 1.4.6]{LW}. Here we only need to note that $R/(x)$
is also a coherent ring (by \cite[Theorem 4.1.1 (1)]{Glaz}) and
over a coherent ring $R$ the $R$-module $\Ext^{n}_{R}(M,N)$ is
coherent for every $n\geq 0$ and every coherent $R$-modules $M$
and $N$ (by \cite[Corollary 2.5.3]{Glaz}).\cqfd\bigskip

Finally, the authors have not been able to extend the third change
of rings theorem for the classical injective dimension
\cite[Theorem 206]{Kap} to the Gorenstein injective dimension over
arbitrary Noetherian rings. However, there are some works which
attempted to give such an extension. See,  for instance,
\cite[Corollary 2.3]{KY} which shows that the desired extension
holds over almost Cohen-Macaulay local rings.
\bigskip

%%%%%%%%%%%%%%%%%%%%%%%%%%%%%%%%%%%%%
\noindent {\bf Acknowledgements.} The authors thank the referee
for his/her careful reading of this work.
 %%%%%%%%%%

\end{section}

%%%%%%%%%%%%%%%%%%%%%%%%%%%%%%%%%%%%%%%%%%%%%%%%%%%%%%%%%%%%
%%%%%%%%%%%%%%%%%%%%%%%%%%%%%%%%%%%%%%%%%%%%%%%%%%%%%%%%%

%%%%%%%%%%%%%%%%%%%%%%%%%%%%%%%%%%%%%%%%%%%%%%%%%%%%%%%%%%%%
%%%%%%%%%%%%%%%%%%%%%%%%%%%%%%%%%%%%%%%%%%%%%%%%%%%%%%%%%
%%%%%%%%%%%%%%%%%%%%%%%%%%%%%%%%%%%%%%%%%%%%%%%%%%%%%%%%%
%%%REFERENCES%%%%%%%%%%%%%%%%%%%%%%%%%%%%%%%%%%%%%%%%%%%%
%%%%%%%%%%%%%%%%%%%%%%%%%%%%%%%%%%%%%%%%%%%%%%%%%%%%%%%%

\bigskip\bigskip

%%%%%%%%%%%%%%%%%%%%%%%%%%%%%%%%%%%%%%%%%%%%%%%%%%%%%%%%
\end{document}